\documentclass[preprint,11pt]{elsarticle}

\makeatletter
\def\ps@pprintTitle{%
\let\@oddhead\@empty
\let\@evenhead\@empty
\def\@oddfoot{}%
\let\@evenfoot\@oddfoot}
\makeatother

\usepackage[titletoc,title]{appendix}
\usepackage{threeparttable}

\makeatletter
  \@addtoreset{chapter}{part}
  \@addtoreset{@ppsaveapp}{part}
\makeatother

\usepackage[latin1]{inputenc}
\usepackage[english]{babel}
\usepackage{caption}

\usepackage{amsmath,amsthm}
\usepackage[nameinlink]{cleveref}
\usepackage{amsfonts}
\usepackage{caption}
\usepackage{relsize}
\usepackage{algorithm}
\usepackage{algorithmicx}
\usepackage{algpseudocode}
\usepackage{xcolor}
\usepackage{cancel}

\usepackage{algcompatible}
\usepackage{setspace}
\usepackage{euscript}
\usepackage{enumerate}
\usepackage[inline]{enumitem}
\usepackage{graphicx}
\usepackage{tikz}
\usepackage{xspace}

\newtheorem{Proposition}{Proposition}

\usetikzlibrary{positioning,chains,fit,shapes,calc}

\usepackage{natbib}
 \bibpunct[, ]{(}{)}{,}{a}{}{,}%

\newtheorem{Property}{Property}
\newtheorem{Remark}{Remark}

\crefname{Proposition}{Proposition}{Propositions}
\crefname{Property}{Property}{Property}
\crefname{algorithm}{Algorithm}{Algorithms}
\crefname{figure}{Figure}{Figures}

\newcommand{\IIO}{\emph{Iterated Inside Out}\xspace}
\newcommand{\NSB}{NS-BPPH\xspace}
\newcommand{\NSBDCS}{NS-BDCS\xspace}

\usetikzlibrary{positioning,chains,fit,shapes,calc}

\begin{document}

\renewcommand{\baselinestretch}{1.0}

\pagestyle{plain}

\begin{frontmatter}

\title{ITERATED INSIDE OUT: \\ a new exact algorithm \\ for the transportation problem}
\date{}

  \author[1]{Roberto Bargetto}
 \author[1,2]{Federico Della Croce}
 \author[1]{Rosario Scatamacchia}
 
  \address[1]{\small Dipartimento di Ingegneria Gestionale e della Produzione, Politecnico di Torino,\\ Corso Duca degli Abruzzi 24, 10129 Torino, Italy, \\{\tt 
\{roberto.bargetto,federico.dellacroce,rosario.scatamacchia\}@polito.it }}
 \address[2]{CNR, IEIIT, Torino, Italy}

\begin{abstract}
We propose a novel exact algorithm for the transportation problem, one of the paradigmatic network optimization problems.
The algorithm, denoted \IIO, requires in input a basic feasible solution and is composed by two main phases that are iteratively repeated until an optimal basic feasible solution is computed. 
In the first ``inside'' phase, the algorithm progressively improves upon a given basic solution by increasing the value of several non-basic variables with negative reduced cost. This phase typically outputs a non-basic feasible solution interior to the constraint set polytope. The second ``out'' phase moves in the opposite direction by iteratively setting to zero several variables until a new improved basic feasible solution is reached.  
Extensive computational tests show that the proposed approach strongly outperforms all versions of network and linear programming algorithms available in the commercial solvers Cplex and Gurobi and other exact algorithms available in the literature. 
\end{abstract}


\end{frontmatter}

\section{Introduction}
\label{Intro}

We consider the transportation problem (TP), one of the historical network optimization problems in the mathematical and operations research communities. In the transportation problem, a given commodity has to be shipped from a number of sources to a number of destinations at minimum cost. The problem can be formalized as follows. Let $ M$ and $ N$ be the set of sources and the set of destinations, respectively, and let $ a_i$ and $ b_j$ denote the level of supply at each source $ i \in M$ and the amount of demand at each destination $ j \in N$. We denote by $c_{ij}$ the unit transportation cost from source $ i \in M$ to destination $ j \in N$. Let $x_{ij} \ge 0$ be a non-negative real variable representing the quantity sent from source $ i \in M$ to destination $ j \in N$.  It is well-known that TP can be always formulated 
such that  $\sum_{i \in M} a_i = \sum_{j \in N} b_j$. Assuming this condition, a linear programming (LP) formulation for TP reads

\begin{align}
\min \quad & z = \sum_{i \in M} \sum_{j \in N} c_{ij} x_{ij} &   \label{o1}\\
& \sum_{j \in N} x_{ij} = a_i  \quad \forall \; i \in M \;\;  \label{c1}\\ 
& \sum_{i \in M} x_{ij} = b_j \quad  \forall \; j \in N \;\;  \label{c2}\\ 
& x_{ij} \ge 0 \quad  \forall \; i \in M, \, j \in N. \; \; \label{c3}
\end{align}

The constraint matrix of model (\ref{o1})--(\ref{c3}) is sparse as each variable appears in just two constraints.  We recall, see, e.g., \cite{Luenberger}, that the constraint matrix is also totally unimodular and hence, if all $a_i$ and $b_j$ are integer, then all basic solutions are integer. Since the best basic feasible solution provides an optimal solution, TP can be solved by any efficient LP solver 
even when the transported quantities are required to be integer ($x_{ij} \in \mathcal{N}$). TP can also be represented as a flow problem on a bipartite graph $G(M,N,E)$, where sources are represented by vertices $i \in M$,
destinations are represented by vertices $j \in N$, and each variable $x_{ij} \geq 0$ corresponds to the flow on edge $(i,j) \in E$ with unit cost $c_{ij}$. 
We also recall, see always \cite{Luenberger}, that, if $\sum_{i \in M} a_i = \sum_{j \in N} b_j$, any basic solution
of model (\ref{o1})--(\ref{c3}) has exactly $|M|+|N|-1$ basic variables (as there is one redundant constraint) and 
corresponds to a spanning tree in $G(M,N,E)$.


TP was stated for the first time by \cite{Monge}. 
Since then, TP has been intensively studied particularly in the twentieth century. We mention here the pioneering works by \cite{Hitchcock}, \cite{Kantorovich}, \cite{Dantzig}, \cite{FordFulkerson}. \cite{Dantzig} provided the first primal simplex algorithm for TP. 
Later, at the end of the twentieth century, several polynomial time algorithms were proposed for the minimum cost flow problem (MCFP), which also generalizes TP, e.g., the primal network simplex algorithm proposed in \cite{Orlin} and the dual network simplex algorithm proposed in \cite{Armstrong}. \cite{Kovacs} indicates that the primal network simplex algorithm is the best performing algorithm for MCFP on dense graphs. \cite{Schrieber} indicate that approaches based on the simplex algorithm (in its various expressions: primal, dual, network) are the best performing approaches for TP.  We finally mention that TP has recently attracted a significant attention in computer vision and machine learning applications, where typically large instances of TP need to be efficiently solved to compute distances between probability measures. We refer to \cite{BaGuaVe} and the references therein for a comprehensive overview on the matter. \cite{Schrieber} introduced a set of benchmark instances of the transportation problem deriving from applications in image processing. These applications call for efficient solution methods for TP, possibly exploiting specific geometrical structures of the transportation costs. 


In this work, we propose a novel exact method for TP. 
The algorithm, denoted \IIO (IIO), requires in input a basic feasible solution and is composed by two main phases that are iteratively repeated until an optimal basic feasible solution is computed. The proposed algorithm exploits the sparseness of the constraint matrix with the practical efficiency of the pivoting operations to explore non-basic solutions, and it limits the iterative computation of the multipliers with respect to the network simplex. The proposed algorithm turns out to be extremely efficient and strongly outperforms all the current state-of-the-art approaches, including all versions of network and linear programming algorithms available in the commercial solvers Cplex and Gurobi.


The remainder of the paper is organized as follows. We briefly recall the network simplex algorithm for TP and introduce the proposed \IIO algorithm in Section \ref{sec:iio}. We discuss further features to speed up the algorithm in Section \ref{sec:speed}. We present extensive computational results on randomly generated instances and benchmark instances from the literature in Section \ref{sec:expe}. Section \ref{sec:concl} concludes the paper with final remarks.

\section{\IIO algorithm}
\label{sec:iio}
\subsection{Preliminaries}
We briefly recall the main steps of the network simplex algorithm for TP, for which we refer to \cite{Luenberger} and its relevant notation, and then introduce our \IIO algorithm.

\paragraph{The network simplex algorithm for TP.}
\label{par:netsimplex}
The algorithm takes as input a feasible basis $B$ and the related vector of non-negative basic variables $\mathbf{x_B}$ corresponding to a spanning tree on the associated graph $G$. The algorithm first computes the simplex multipliers related to $B$: the multipliers $u_i$ for each source $i \in M$ and the multipliers $v_j$ for each destination $j \in N$. Since each multiplier can be computed in constant time and all $|M|+|N|$ multipliers have to be computed, this step requires $\Theta(|M|+|N|)$ time complexity. Then, the algorithm computes the reduced cost $r_{ij} = c_{ij} - u_i - v_j$ of each non-basic variable $x_{ij} \not \in \mathbf{x_B}$ in constant time. If all reduced costs are non-negative, the basis $B$ is optimal and the algorithm terminates. A worst-case $O(|M||N|)$ time complexity holds for this step as, potentially, all reduced costs need to be computed. 

However, the average complexity per iteration may be much less since, whenever a negative reduced cost of a non-basic variable is found, the algorithm may proceed to the following pivoting step. In the pivoting step, a basic variable is replaced by a non-basic variable $x_{ij}$ with negative reduced cost by searching for the unique cycle in graph $G$ induced by the spanning tree associated with $\mathbf{x_B}$ and edge $(i,j)$. This task corresponds to computing the unique simple path $P$ between source $i$ and destination $j$ in the spanning tree. 


The time complexity of a pivoting operation is $\Theta(|P|)$, where $|P|$ denotes the size of the path (i.e., the number of nodes included in it) and $|P| < |M|+|N|$. Typically, $|P|$ is much smaller than $|M|+|N|$ so that, in practice, the computational effort of a pivoting operation is negligible with respect to that of the multipliers computation step. After a pivoting operation, a new basic solution is obtained and the algorithm iterates. It turns out that the time complexity of an iteration of the network simplex algorithm for TP is lower bounded by $\Omega(|M|+|N|)$, namely the complexity of the multipliers computation step. 

The following property trivially holds for TP.

\begin{Property}
\label{Property1} \hfill

\noindent
Consider an instance of TP with supplies $a_i$ ($i \in M$) and demands $b_j$ ($j \in N$). 
\begin{itemize}
\item[a)] Given a non-basic feasible solution $\mathbf{x}$, 
any subset of variables $\mathbf{\mathbf{x_B}} \subset \mathbf{x}$ with cardinality $|M|+|N|-1$, such that the corresponding edges in graph $G$ form a spanning tree, 
constitutes a basic feasible solution of a modified instance of TP with the following supplies and demands

\begin{center}
$a'_i=a_i - \sum\limits_{j \in N: \; x_{ij}\not \in \mathbf{x_B}}x_{ij} \;\;\; \forall i \in M$,\\
 $b'_j=b_j - \sum\limits_{i \in M: \; x_{ij}\not \in \mathbf{x_B}}x_{ij} \;\;\; \forall j \in N$.
\end{center}

\item[b)] Consider a basic feasible solution $\mathbf{x_B}$. 
Let modify the instance by increasing $a_h$ and $b_l$ by a given value $\gamma > 0$ so that  $a_h = a_h +\gamma$ and $b_l = b_l +\gamma$ where $x_{hl} \not \in \mathbf{x_B}$. Correspondingly, generate a non-basic feasible solution to this new instance by setting  $x_{hl} = \gamma$.  
Then, it is always possible to restore a basic feasible solution of the modified instance by applying a pivoting operation on the unique cycle induced in graph $G$ by the spanning tree associated with $\mathbf{x_B}$ plus edge $(h,l)$.
\end{itemize}
\end{Property}

\subsection{The proposed approach}
The proposed \IIO algorithm exploits \cref{Property1} together with the sparseness of the constraint matrix and the practical efficiency of a pivoting operation. Given a basic feasible solution $\mathbf{x_B}$, we recall that the objective function $z$ can be expressed as a function of the non-basic variables and their reduced costs, namely $z = z_0 + \sum_{i,j: x_{ij} \not \in \mathbf{x_B}}r_{ij}x_{ij}$, where $z_0$ indicates the transportation costs related to the basic variables. Notice that, if the non-basic variables with $r_{ij} \geq 0$ are not considered, the more the non-basic variables with $r_{ij} < 0$ increase their value, the more the new solution improves its objective. 


A first attempt in this direction, taking into account the sparseness of the constraint matrix, was proposed in \cite{Bulut}, where several disjoint cycles involving distinct non-basic variables with negative reduced costs are computed within a single iteration of the so-called multi-loop simplex algorithm. That approach, however, was capable of reducing only marginally the computational time with respect to a standard implementation of the simplex algorithm. 


We propose a different approach where, starting from a basic solution $\mathbf{x_B}$, we compute, with a limited computational effort requiring only pivoting operations, an improved feasible non-basic solution $\mathbf{\hat{x}}$. In that solution, compared to $\mathbf{x_B}$, many originally non-basic variables with $r_{ij} < 0$ are now strictly positive while all the other non-basic variables are equal to $0$. From the non-basic solution $\mathbf{\hat{x}}$, 
we reach a new basic feasible solution $\mathbf{x_{B'}}$ 
again through repeated pivoting operations. The process requires the computation of reduced costs and multipliers for solutions $\mathbf{x_B}$  and $\mathbf{x_{B'}}$, but avoids that computation for all intermediate solutions computed between $\mathbf{x_B}$ and $\mathbf{x_{B'}}$. 


\IIO receives as input an initial basic feasible solution and is split into two main phases that are iteratively repeated until an optimal basic feasible solution is computed. The two phases of the algorithm are described below.

\medskip

\paragraph{Phase 1: the ``inside'' phase.} In the first phase, the current basic feasible solution is progressively improved by considering one at a time the non-basic variables with negative reduced cost. For each  non-basic variable $x_{ij}$, the related value is increased as much as possible to some given amount $k$ by means of a pivoting operation where a basic variable is set to 0 but it is not removed from the basic solution. This operation, see \cref{Property1} (a), corresponds to considering a modified problem where we fictitiously remove the variable $x_{ij}$: correspondingly, we decrease $a_i$ and $b_j$ by the amount $k$  ($a_i = a_i-k$, $b_j = b_j-k$) so as to refer always to the same basis. Also, each pivoting operation requires the search of a path $P$  in the same spanning tree in $G$ with time complexity $\Theta(|P|)$ with no need to recompute the simplex multipliers and the reduced costs.


 We remark that, during the pivoting operations, the values of the initial basic variables progressively vary and many of them typically oscillate between value $0$ (that can be reached more than once) and various positive values. At the end of this phase, we obtain a feasible solution where, typically, the number of variables with strictly positive value is typically greater than the size of the basis $|M|+|N|-1$, namely a  non-basic feasible solution that is {\em inside} the feasibility region determined by the constraint set. 


\paragraph{Phase 2: the ``out'' phase.}

In the second phase, we consider the variables added in Phase 1 with the aim of leaving the interior part of the polytope and getting back to a basic feasible solution. The phase starts with the initial basis considered in Phase 1 where the basic variables keep the value reached at the end of that phase. 
Then, one at a time, all the added non-basic variables are considered exploiting \cref{Property1} (b). We remark that these variables are strictly positive. 
Consider the first non-basic variable $x_{hl}$ added in Phase 1. In order to perform a pivoting operation, we identify the related path $P$ in the spanning tree in $G$ corresponding to the current basis. Since $r_{hl} < 0$, it would be indeed convenient to increase the value of $x_{hl}$ in order to further improve the objective function. If variable $x_{hl}$ can increase, the basis changes with one of the basic variables leaving the basis and $x_{hl}$ entering the basis. If $x_{hl}$ cannot increase in case of degeneracy, $x_{hl}$ would keep the same value but still would enter the basis. After executing the first pivoting operation, we get a new basic feasible solution to the problem with updated values of $a_h$ and $b_l$. From the second pivoting operation on, the reduced costs of the other variables in Phase 1 are no more meaningful. Hence, we iteratively evaluate whether it is convenient to increase or to decrease the value of the next non-basic variable $x_{pq}$ by considering the cost coefficients of the variables involved in the corresponding cycle (e.g., suppose we have a cycle with variables $x_{pq}$, $x_{pr}$, $x_{sr}$, $x_{sq}$. If $c_{pq} - c_{pr} + c_{sr} - c_{sq} < 0$, then it is convenient to increase $x_{pq}$. Else, if $c_{pq} - c_{pr} + c_{sr} - c_{sq} > 0$, it is convenient to decrease $x_{pq}$. If $c_{pq} - c_{pr} + c_{sr} - c_{sq} = 0$, no improvement can be obtained and we arbitrarily try to increase  $x_{pq}$ to reach a new basis).
If $x_{pq}$ should increase, the same analysis of variable $x_{hl}$ applies. Else, either $x_{pq}$ is decreased to value $0$, while the basic variables remain the same but change their value, or $x_{pq}$ is decreased until a basic variable reaches value 0 and leaves the basis. In the latter case, again, $x_{pq}$  enters the basis. Correspondingly, a new basic feasible solution is obtained to the problem with updated values of $a_p$ and $b_q$. The approach is iterated until all the variables added in Phase 1 have been taken into account.
At the end of Phase 2, a (generally improved) basic feasible solution of the original problem is obtained. \\

The basic solution obtained at the end of Phase 2 requires then the computation of the corresponding multipliers and reduced costs. If all reduced costs are non-negative, the basic solution is optimal and the algorithm terminates. Otherwise, the algorithm re-applies the two phases.

\begin{Remark}
In case of degenerate basic solutions, 
it may occur that the solution value of the basic solution at the beginning of Phase 1 is not improved after applying 
all pivoting steps of the two phases. In this case, \IIO may eventually cycle (like the simplex algorithm) among basic solutions without converging to an optimal solution. 
While in our computational testing such behavior never occurred, it is sufficient to add 
an anti-cycling mechanism, such as the strongly feasible bases in \cite{Cunningham}, whenever Phase 1 and Phase 2 do not yield an improved basic solution. 
Once an improvement is found, the algorithm applies the two phases again.
\end{Remark}

The pseudo-code of the proposed approach (where the anti-cycling mechanism is omitted) is depicted in the following \cref{IIO}.
\begin{algorithm}[ht!]
\caption{\IIO}
\label{IIO}
\begin{algorithmic}[1]
\State{\textbf{Input:} a basic feasible solution $\mathbf{x_B}$.}
\State 
Compute for each source $a_i$ the related multiplier $u_i$; compute for each destination $b_j$ the related multiplier $v_j$.
\State 
Compute the reduced costs $r_{ij} = c_{ij} - u_i - v_j$ of the non-basic variables. 
\State \textbf{If} $r_{ij} \geq 0 \; \; \forall x_{ij} \not \in \mathbf{x_B}$ \textbf{then} $\mathbf{x_B}$ is optimal, \textbf{return} $\mathbf{x_B}$.\\
\textbf{Else} \\
\emph{Phase 1:} Increase (when possible) one at a time the value of the non-basic variables with reduced cost $<0$ by pivoting operations. 
\State 
\emph{Phase 2:} Given the basic solution considered in Phase 1, analyze one at a time the added variables and apply pivoting operations leading to a new feasible basic solution  $\mathbf{x_B}$. 
\State 
\textbf{Go to} 2. 
\State \textbf{End If}
\end{algorithmic}
\end{algorithm}

\begin{Proposition}
\label{Prop0}
\IIO solves TP to optimality.
\end{Proposition}

\begin{proof}
\IIO requires in input a basic feasible solution as a standard network simplex algorithm. Then, at the beginning of Phase 1, the first pivoting step corresponds to a pivoting step of the network simplex algorithm. Hence, unless degeneracy occurs,  
the algorithm computes an improved basic solution before re-applying Phase 1. 
Correspondingly, in absence of degeneracy, the repeated application of the two phases guarantees the convergence of the algorithm to an optimal solution. In case of degeneracy, the addition of any anti-cycling mechanism allows \IIO to obtain an improved basic solution and correspondingly to reach an optimal solution in the next iterations.
\end{proof}
\subsection{An illustrative example}
\label{exIIO}

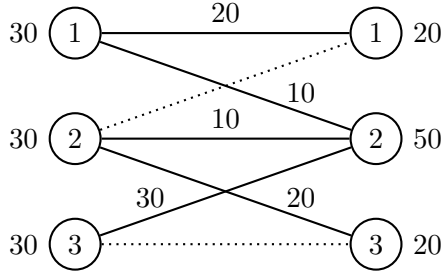
\begin{figure}[!h]
\centering
\begin{tikzpicture}[thick,
  fsnode/.style={fill=white, draw,circle},
  ssnode/.style={fill=white, draw,circle},
]

\begin{scope}[start chain=going below,node distance=7mm]

  \node[fsnode,on chain] (f1) [label=left: $30$] {1};
   \node[fsnode,on chain] (f2) [label=left: $30$] {2};
 \node[fsnode,on chain] (f3) [label=left: $30$] {$3$};
\end{scope}

\begin{scope}[xshift=4cm,start chain=going below,node distance=7mm]

  \node[ssnode,on chain] (s1) [label=right: $20$] {1};
  \node[ssnode,on chain] (s2) [label=right: $50$] {2};
   \node[ssnode,on chain] (s3) [label=right: $20$] {$3$};
\end{scope}

\draw[] (f1) -- (s1)  node[pos=0.5, above] {$20$ };
\draw (f1) -- (s2)  node[pos=0.8, above] {$10$ };
{\draw[dotted] (f2) -- (s1)node[pos=0.2, above]{} ;}
\draw (f2) -- (s2)  node[pos=0.5, above] {$10$ };
\draw (f2) -- (s3)  node[pos=0.8, above] {$20$};
\draw (f3) -- (s2)  node[pos=0.2, above] {$30$};
{\draw[dotted] (f3) -- (s3)node[pos=0.2, above]{} ;}
\end{tikzpicture}
\caption{ The starting basic solution on graph $G$.}
\label{Exfig1}
\end{figure}

Consider an  instance of TP where $M = N = 3$ with sources supplies $a = [30 \; 30 \; 30]$
and destinations demands $b = [20 \; 50 \; 20]$. The transportation costs are: $c_{11} = 5, c_{12} = 1, c_{13} = 7, c_{21} = 1, c_{22} = 1, c_{23} = 5, c_{31} = 6, c_{32} = 1, c_{33} = 2$. Let the initial basic solution $\mathbf{x_B}$ be constituted by variables $x_{11}=20$, $x_{12}=10$, $x_{22}=10$, $x_{23}=20$, $x_{32}=30$ with objective function $z=250$. Concerning the non-basic variables, $x_{21}$ and $x_{33}$ have negative reduced costs ($r_{21} = -4$ and $r_{33} =-3$, respectively) while $x_{13}$ and $x_{31}$ have positive reduced costs ($r_{13} = 2$ and $r_{31} =1$, respectively). In \cref{Exfig1} is represented the solution on the corresponding bipartite graph where edges related to variables with positive reduced costs are omitted and edges related to variables with negative reduced costs are dotted.

\begin{figure}[!h]
\centering
\begin{tikzpicture}[thick,
  fsnode/.style={fill=white, draw,circle},
  ssnode/.style={fill=white, draw,circle},
]

\begin{scope}[start chain=going below,node distance=7mm]

  \node[fsnode,on chain] (f1) [label=left: $30$] {1};
   \node[fsnode,on chain] (f2) [label=left: {$20$} $\cancel{30}$] {2};
 \node[fsnode,on chain] (f3) [label=left: $30$] {$3$};
\end{scope}

\begin{scope}[xshift=4cm,start chain=going below,node distance=7mm]

  \node[ssnode,on chain] (s1) [label=right: $\cancel{20} $ {$10$}] {1};
  \node[ssnode,on chain] (s2) [label=right: $50$] {2};
   \node[ssnode,on chain] (s3) [label=right: $20$] {$3$};
\end{scope}



\draw[] (f1) -- (s1) node[pos=0.5, above] {$10$ };
\draw (f1) -- (s2) node[pos=0.8, above] {$20$ };
\draw (f2) -- (s2)node[pos=0.5, above] {$0$ };
\draw (f2) -- (s3)  node[pos=0.8, above] {$20$};
\draw (f3) -- (s2)  node[pos=0.2, above] {$30$};
{\draw[dotted] (f3) -- (s3);}
\end{tikzpicture}
\caption{ Phase 1 after the first pivoting operation.}
\label{Exfig2}
\end{figure}
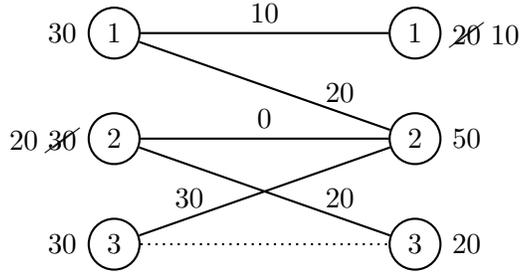

As indicated in the algorithm pseudo-code, in Phase 1 we consider, one at a time, the non-basic variables with negative reduced cost; that is, in this case, first $x_{21}$ and then $x_{33}$. Like in a standard pivoting operation, we increase as much as possible the value of $x_{21}$ until a basic variable reaches value $0$. Since we have a cycle involving the variables $x_{21}$, $x_{22}$, $x_{12}$, $x_{11}$, where the basic variables $x_{11}$ and $x_{22}$ decrease while the basic variable $x_{12}$ increases, we obtain $x_{21}=10$. Then, exploiting \cref{Property1} (a), we fictitiously remove variable $x_{21}$ from the problem by correspondingly updating $a_2$ and $b_1$. As shown in \cref{Exfig2}, the basic variables remain the same with updated values $x_{11}=10$, $x_{12}=20$, $x_{22}=0$, $x_{23}=20$, $x_{32}=30$ and $z=210$.

\begin{figure}[!h]
\centering
\begin{tikzpicture}[thick,
  fsnode/.style={fill=white, draw,circle},
  ssnode/.style={fill=white, draw,circle},
]

\begin{scope}[start chain=going below,node distance=7mm]

  \node[fsnode,on chain] (f1) [label=left: $30$] {1};
   \node[fsnode,on chain] (f2) [label=left: {$20$} $\cancel{30}$] {2};
   \node[fsnode,on chain] (f3) [label=left: {$10$} $\cancel{30}$] {3};
\end{scope}

\begin{scope}[xshift=4cm,start chain=going below,node distance=7mm]

  \node[ssnode,on chain] (s1) [label=right: $\cancel{20} $ {$10$}] {1};
  \node[ssnode,on chain] (s2) [label=right: $50$] {2};
  \node[ssnode,on chain] (s3) [label=right: $\cancel{20} $ {$0$}] {3};
\end{scope}



\draw[] (f1) -- (s1) node[pos=0.5, above] {$10$ };
\draw (f1) -- (s2) node[pos=0.8, above] {$20$ };
\draw (f2) -- (s2)node[pos=0.5, above] {$20$ };
\draw (f2) -- (s3)  node[pos=0.8, above] {$0$};
\draw (f3) -- (s2)  node[pos=0.2, above] {$10$};
\end{tikzpicture}
\caption{ Phase 1 after the second pivoting operation.}
\label{Exfig3}
\end{figure}
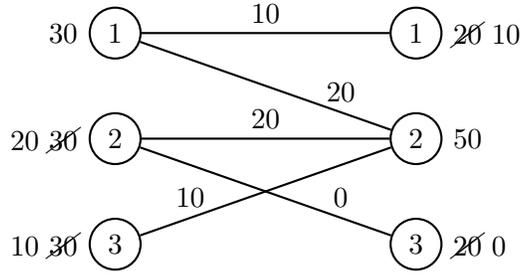

We now consider the variable $x_{33}$. Analyzing the cycle with variables $x_{33}$, $x_{32}$, $x_{22}$, $x_{23}$, we obtain $x_{33}=20$. We fictitiously remove variable $x_{33}$ from the problem by updating $a_3$ and $b_3$. As indicated in \cref{Exfig3}, the basic variables remain the same with updated values $x_{11}=10$, $x_{12}=20$, $x_{22}=20$, $x_{23}=0$, $x_{32}=10$ and $z=150$.

Phase 1 terminates with a non-basic feasible solution where the number of strictly positive variables is $6$ which is greater than the size of $\mathbf{x_B}$ ($|M|+|N|-1=5$). We have $x_{11}=10$, $x_{12}=20$, $x_{21}=10$, $x_{22}=20$, $x_{32}=10$,
$x_{33}=20$ and $z=150$.

\begin{figure}[!h]
\centering
\begin{tikzpicture}[thick,
  fsnode/.style={fill=white, draw,circle},
  ssnode/.style={fill=white, draw,circle},
]

\begin{scope}[start chain=going below,node distance=7mm]

  \node[fsnode,on chain] (f1) [label=left: $30$] {1};
   \node[fsnode,on chain] (f2) [label=left: $30$] {2};
   \node[fsnode,on chain] (f3) [label=left: {$10$} $\cancel{30}$] {3};
\end{scope}

\begin{scope}[xshift=4cm,start chain=going below,node distance=7mm]

  \node[ssnode,on chain] (s1) [label=right: $20$] {1};
  \node[ssnode,on chain] (s2) [label=right: $50$] {2};
  \node[ssnode,on chain] (s3) [label=right: $\cancel{20} $ {$0$}] {3};
\end{scope}



\draw[] (f1) -- (s1)  node[pos=0.5, above] {$10$};
\draw (f1) -- (s2)  node[pos=0.8, above] {$20$};
{\draw[dotted] (f2) -- (s1)node[pos=0.2, above]{$10$} ;}
\draw (f2) -- (s2)  node[pos=0.5, above] {$20$};
\draw (f2) -- (s3)  node[pos=0.8, above] {$0$};
\draw (f3) -- (s2)  node[pos=0.2, above] {$10$};
\end{tikzpicture}
\caption{Phase 2 before processing the first non-basic variable $x_{21}$.}
\label{Exfig5}
\end{figure}
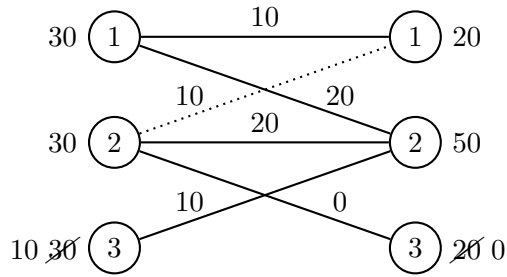

At the beginning of Phase 2, we consider the initial basic solution $\mathbf{x_B}$ with variables 
$x_{11}$, $x_{12}$, $x_{22}$, $x_{23}$, $x_{32}$ (having their last updated values previously indicated) and the first of the added non-basic variables with a strictly positive value, namely $x_{21}$. Exploiting \cref{Property1} (b), we modify the current values of supply $a_2$ and demand $b_1$, namely $a_2 = 20 + x_{21} = 30$, $b_1 = 10 + x_{21} = 20$ (see \cref{Exfig5} where the edge related to $x_{21}$ is dotted). 

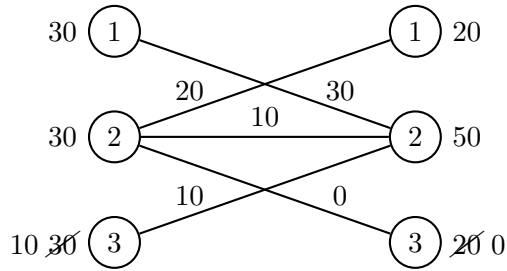
\begin{figure}[H]
\centering
\begin{tikzpicture}[thick,
  fsnode/.style={fill=white, draw,circle},
  ssnode/.style={fill=white, draw,circle},
]

\begin{scope}[start chain=going below,node distance=7mm]

  \node[fsnode,on chain] (f1) [label=left: $30$] {1};
   \node[fsnode,on chain] (f2) [label=left: $30$] {2};
   \node[fsnode,on chain] (f3) [label=left: {$10$} $\cancel{30}$] {3};
\end{scope}

\begin{scope}[xshift=4cm,start chain=going below,node distance=7mm]

  \node[ssnode,on chain] (s1) [label=right: $20$] {1};
  \node[ssnode,on chain] (s2) [label=right: $50$] {2};
  \node[ssnode,on chain] (s3) [label=right: $\cancel{20} $ {$0$}] {3};
\end{scope}



\draw (f1) -- (s2)  node[pos=0.8, above] {$30$ };
{\draw (f2) -- (s1)node[pos=0.2, above]{$20$} ;}
\draw (f2) -- (s2)  node[pos=0.5, above] {$10$};
\draw (f2) -- (s3)  node[pos=0.8, above] {$0$};
\draw (f3) -- (s2)  node[pos=0.2, above] {$10$};
\end{tikzpicture}
\caption{ Phase 2 after processing the first non-basic variable $x_{21}$.}
\label{Exfig6}
\end{figure}

Variable $x_{21}$ induces the corresponding cycle with variables $x_{21}$, $x_{22}$, $x_{12}$, $x_{11}$. It is convenient to increase the value of variable $x_{21}$ as $r_{21} = -4$. The new values of the variables in the cycle become $x_{12}=  30$, $x_{21}=  $20$, x_{22} =  10$ and $x_{11}=0$ with $z=110$. Hence, variable $x_{11}$ leaves the basis and the new basic variables are 
$x_{12}$, $x_{21}$, $x_{22}$, $x_{23}$ and $x_{32}$, as indicated in \cref{Exfig6}.

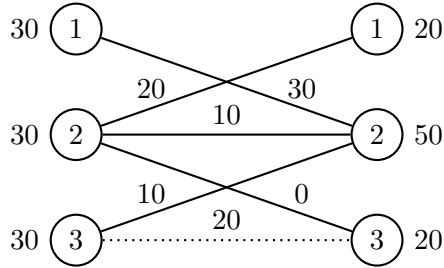
\begin{figure}[H]
\centering
\begin{tikzpicture}[thick,
  fsnode/.style={fill=white, draw,circle},
  ssnode/.style={fill=white, draw,circle},
]

\begin{scope}[start chain=going below,node distance=7mm]

  \node[fsnode,on chain] (f1) [label=left: $30$] {1};
   \node[fsnode,on chain] (f2) [label=left: $30$] {2};
 \node[fsnode,on chain] (f3) [label=left: $30$] {$3$};
\end{scope}

\begin{scope}[xshift=4cm,start chain=going below,node distance=7mm]

  \node[ssnode,on chain] (s1) [label=right: $20$] {1};
  \node[ssnode,on chain] (s2) [label=right: $50$] {2};
   \node[ssnode,on chain] (s3) [label=right: $20$] {$3$};
\end{scope}



\draw (f1) -- (s2)  node[pos=0.8, above] {$30$ };
{\draw (f2) -- (s1)node[pos=0.2, above]{$20$} ;}
\draw (f2) -- (s2)  node[pos=0.5, above] {$10$};
\draw (f2) -- (s3)  node[pos=0.8, above] {$0$};
\draw (f3) -- (s2)  node[pos=0.2, above] {$10$};
{\draw[dotted] (f3) -- (s3)node[pos=0.5, above]{$20$} ;}
\end{tikzpicture}
\caption{ Phase 2 before processing the second non-basic variable $x_{33}$.}
\label{Exfig6b}
\end{figure}

Then, we consider the added variable $x_{33}$ restoring the original values of supply $a_3$ and demand $b_3$, as shown in \cref{Exfig6b} where the edge related to $x_{33}$ is dotted.

\begin{figure}[!h]
\centering
\begin{tikzpicture}[thick,
  fsnode/.style={fill=white, draw,circle},
  ssnode/.style={fill=white, draw,circle},
]

\begin{scope}[start chain=going below,node distance=7mm]
+
  \node[fsnode,on chain] (f1) [label=left: $30$] {1};
   \node[fsnode,on chain] (f2) [label=left: $30$] {2};
 \node[fsnode,on chain] (f3) [label=left: $30$] {$3$};
\end{scope}

\begin{scope}[xshift=4cm,start chain=going below,node distance=7mm]

  \node[ssnode,on chain] (s1) [label=right: $20$] {1};
  \node[ssnode,on chain] (s2) [label=right: $50$] {2};
   \node[ssnode,on chain] (s3) [label=right: $20$] {$3$};
\end{scope}



\draw (f1) -- (s2)  node[pos=0.8, above] {$30$ };
{\draw (f2) -- (s1)node[pos=0.2, above]{$20$} ;}
{\draw (f2) -- (s2)  node[pos=0.5, above] {$10$ };}
\draw (f3) -- (s2)  node[pos=0.2, above] {$10$ };
{\draw (f3) -- (s3)node[pos=0.5, above]{$20$ } ;}
\end{tikzpicture}
\caption{ Phase 2 after processing the second non-basic variable $x_{33}$.}
\label{Exfig7}
\end{figure}
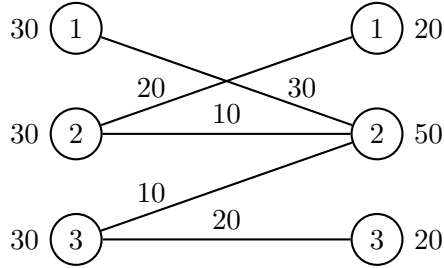

The presence of variable $x_{33}$ induces 
the cycle with variables $x_{33}$, $x_{32}$, $x_{22}$, $x_{23}$. 
Analyzing the cycle, the added variable $x_{33}$ should increase as $c_{33} - c_{32} + c_{22} - c_{23} = -3 < 0$.  But since $x_{23}=0$, $x_{33}$ cannot increase; hence, the solution remains the same with $z=110$, where the variable $x_{33}$ enters the basis and the variable $x_{23}$ leaves the basis, as shown in \cref{Exfig7}.

Phase 2 terminates with a new feasible basic solution $\mathbf{x_B}$ composed by variables $x_{12}=  30$, $x_{21}=  20$, $x_{22}=  10$, $x_{32}=10$, $x_{33}=20$ and $z=110$. At this point, as all reduced costs are non-negative, the algorithm stops having reached an optimal solution.

\section{Speeding up the algorithm}
\label{sec:speed}
Several issues have to be considered in the application of \IIO. First, an initial basic feasible solution needs to be determined with the corresponding computation of the multipliers and the reduced costs of the non-basic variables. While the computation of all multipliers is unavoidable,  we aim to limit the number of reduced costs to be considered during the iterations (apart from the last iteration), since computing all the reduced costs requires $\Theta(|M||N|)$ time complexity, and this complexity may slow down the performance of the algorithm. 


Besides, several basic variables progressively vary in Phase 1 and the more non-basic variables are set to a positive value, the larger the number of basic variables that progressively reach value $0$. If a basic variable with value 0 is supposed to decrease its value during a subsequent pivoting operation, this decrease cannot occur, the considered non-basic variable cannot reach a value greater than zero, and the solution does not change. Thus, it is desirable to avoid, when possible, the computation of cycles that cannot lead to different solutions in the first phase of the algorithm. We discuss all the above mentioned issues in the following subsections.

\subsection{Computing an initial feasible solution}
We tested four well-known approaches from the literature to generate a basic feasible solution, namely 
the North-West Corner (NWC) method, Vogel's Approximation (VA) method, the Matrix Minimum Rule (MMR) method and
the Tree Minimum Rule (TMR) method. We refer to \cite{Schwinn} for a recent survey on heuristic methods for the transportation problem. In agreement with \cite{Schwinn}, the best compromise between solution quality and computational time was provided by MMR in our experimental tests. 

\subsection{Limiting the number of reduced costs to be computed in Phase 1}
This matter has already been tackled in other publications on TP focusing on the simplex algorithm, see, e.g., \cite{Gottschlich} and \cite{Schmitzer}, where just one negative reduced cost is necessary to process one iteration. In particular, in the so-called shortlist method presented in \cite{Gottschlich}, only a subset (shortlist) of the variables with least transportation costs is considered for each source. A set of candidate variables with negative reduced cost is selected within the shortlists, and the variable with the most negative reduced cost is chosen for a pivoting operation. If no improvement can be achieved within the shortlists, all variables are considered and the method performs an iterative analysis of the sources to select the next variables for pivoting operations. We refer the reader to \cite{Gottschlich} for more details on the shortlist method.


We apply a simplified version of the shortlist method where we just introduce a parameter $\alpha$, which represents the size of the subset of the variables with the least transport costs to be taken into account in the application of \IIO. In the same spirit of \cite{Gottschlich}, once all non-basic variables in the subset have non-negative reduced cost, also all other variables are considered in the next iterations of \IIO.

\subsection{Disregarding cycles computation involving degenerate basic variables}\label{ss:coloredspat}
In the first phase of \IIO, we are interested in applying a pivoting operation only if a non-basic variable $ x_{ij}$ with negative reduced cost can enter the basis with a strictly positive value, so as to decrease the objective function. This situation occurs either if all basic variables, corresponding to the edges of the unique path from source $i$ to destination $j$ in the spanning tree, have a strictly positive value, or if each basic variable with value $0$ has the related edge in an even position in the path from $i$ to $j$. In the latter case, variable $x_{ij}$ would certainly assume a strictly positive value with a pivoting operation and all basic variables with value $0$ would increase their value as well. 

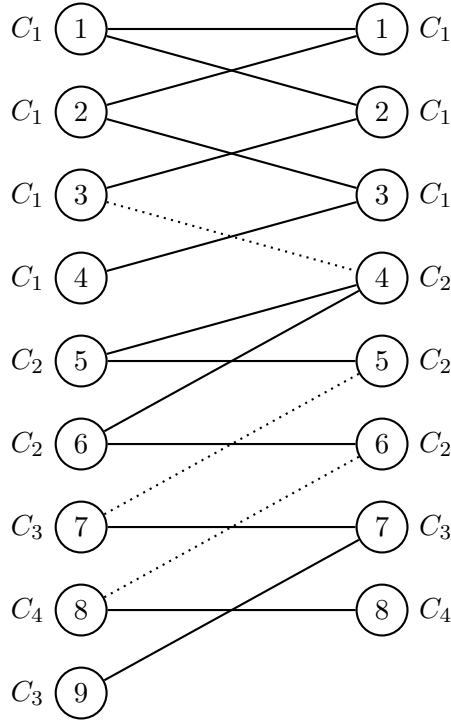
\begin{figure}[ht!]
\centering
\begin{tikzpicture}[thick,
  fsnode/.style={fill=white, draw,circle},
  ssnode/.style={fill=white, draw,circle},
]

\begin{scope}[start chain=going below,node distance=4mm]

  \node[fsnode,on chain] (f1) [label=left: $C_1$]{1};
  \node[fsnode,on chain] (f2) [label=left: $C_1$]{2};
 \node[fsnode,on chain] (f3) [label=left: $C_1$]{3};
 \node[fsnode,on chain] (f4) [label=left: $C_1$]{4};
 \node[fsnode,on chain] (f5) [label=left: $C_2$]{5};
 \node[fsnode,on chain] (f6) [label=left: $C_2$]{6};
 \node[fsnode,on chain] (f7) [label=left: $C_3$]{7};
 \node[fsnode,on chain] (f8) [label=left: $C_4$]{8};
 \node[fsnode,on chain] (f9) [label=left: $C_3$]{9};
\end{scope}

\begin{scope}[xshift=4cm,start chain=going below,node distance=4mm]

  \node[ssnode,on chain] (s1) [label=right: $C_1$]{1};
  \node[ssnode,on chain] (s2) [label=right: $C_1$]{2};
   \node[ssnode,on chain] (s3) [label=right: $C_1$]{3};
   \node[ssnode,on chain] (s4) [label=right: $C_2$]{4};
   \node[ssnode,on chain] (s5) [label=right: $C_2$]{5};
   \node[ssnode,on chain] (s6) [label=right: $C_2$]{6};
   \node[ssnode,on chain] (s7) [label=right: $C_3$]{7};
   \node[ssnode,on chain] (s8) [label=right: $C_4$]{8};
\end{scope}



\draw[] (f1) -- (s1)  node[pos=0.5, above] {};
\draw (f1) -- (s2)  node[pos=0.8, above] {};
\draw (f2) -- (s3)  node[pos=0.8, above] {};
\draw (f2) -- (s1);
\draw (f3) -- (s2);
\draw[dotted] (f3) -- (s4)  node[pos=0.2, above] {};
\draw (f4) -- (s3);
\draw (f5) -- (s4);
\draw (f5) -- (s5);
\draw (f6) -- (s4);
\draw (f6) -- (s6);
\draw (f7)[dotted] -- (s5);
\draw (f7) -- (s7);
\draw (f8)[dotted] -- (s6);
\draw (f8) -- (s8);
\draw (f9) -- (s7);

\end{tikzpicture}
\caption{ An example of colored spanning tree.}
\label{fig:example}
\end{figure}

We handle this issue by means of a coloring of the spanning tree.
Apart from the first iteration where all vertices have the same color and 
a related pivoting operation is applied, 
for all other iterations we color the spanning tree as follows. Any edge $(i,j)$ of the spanning tree corresponds to a basic variable with source $i$ and destination $j$. We say that an edge is a {\em degenerate basic edge} if the related basic variable has value $0$. 
Assuming, w.l.o.g., that source $1$ is the root node of the spanning tree, we assign color $C_1$ to all nodes (sources/destinations) reachable by source $1$ without crossing degenerate basic edges.
Let denote by $T_1$ the subtree induced by all nodes with color $C_1$. Then, following the increasing orders of sources and destinations in the bipartite graph, let $(i,j)$ be the first degenerate edge connecting a source $i$ (destination $j$) with color $C_1$ to a destination $j$ (source $i$) not yet colored. 
We color the destination $j$ (source $i$) with the color $C_2$ and denote $j$ ($i$) as the root node of subtree $T_2$ with color $C_2$. 
Also, we assign the color $C_2$ to all (uncolored) nodes reachable by destination $j$ (source $i$)  without crossing further degenerate basic edges, and denote $T_1$ as the parent of $T_2$. 
With the same approach we color all other nodes by adding colors when necessary, detecting subtrees and the relationships parent/child as indicated above. 
Notice that the coloring of the tree is not performed from scratch after each pivoting operation but it is progressively updated when an edge becomes degenerate. 
Similarly, whenever a degenerate edge connecting $T_i$ and $ T_j$ with $T_i$ parent of $T_j$ becomes non-degenerate, then $T_i$ and $T_j$ are merged and all nodes in $T_j$ 
are assigned color $C_i$.
We provide an example of a colored spanning tree in \cref{fig:example}. 

In the figure, there are $9$ sources and $8$ destinations and the spanning tree corresponding to the current basic solution has three degenerate edges $(3,4)$, $(7,5)$ and $(8,6)$ (dotted in the picture). Correspondingly, 
\begin{enumerate*}[label=(\roman*)]
\item sources 1, 2, 3, 4 and destinations 1, 2, 3 form the first subtree $T_1$ with color $C_1$ and source 1 as root node; 
\item sources 5, 6 and destinations 4, 5, 6 form the second subtree $T_2$ with color $C_2$ and destination 4 as root node; 
\item sources 7, 9 and destination 7 form the third subtree $T_3$ with color $C_3$ and source $7$ as root node; 
\item source 8 and destination 8 form the fourth subtree $T_4$ with color $C_4$ and source $8$ as root node. 
\end{enumerate*}

We remark that, in Phase 1, whenever a degenerate basic edge is generated in the spanning tree, a new subtree is created and all its nodes are assigned a new color. 
While the worst case complexity of this step is $O(|M|+|N|)$, in practice the required computational effort is much smaller. The same consideration holds whenever in Phase $1$ a basic variable with value 0 increases its value and the corresponding edge is no more degenerate. In this case, all nodes of the related subtree merge with the nodes of the parent subtree and take the corresponding color.


Given a colored spanning tree, we can possibly evaluate in constant time whether a given non-basic variable $x_{ij}$ with negative reduced cost would take a value $> 0$ after a pivoting operation, thus avoiding the computation of the related cycle if the variable cannot increase. The following exhaustive cases hold for variable $x_{ij}$:
\begin{enumerate}
	\item\label{samec} Source $i$ and destination $j$ have the same color $C_l$ (they both belong to subtree $T_l$);
	\item\label{diffcpc} Source $i$ has color $C_l$, destination $j$ has color $C_m \neq C_l$ and subtree $T_l$ is either the parent or the child of $T_m$;
	\item\label{diffcbro} Source $i$ has color $C_l$, destination $j$ has color $C_m \neq C_l$ and the subtrees $T_l$ and $T_m$ have the same parent, namely they are siblings.
	\item\label{diffcboh} Source $i$ and destination $j$ have different colors and none of the above cases \ref{diffcpc}, \ref{diffcbro} holds.
\end{enumerate}

The following proposition holds.
\begin{Proposition}
\label{Prop1}
Given a colored spanning tree and a non-basic variable $x_{ij}$, if one of the above cases \ref{samec},\ref{diffcpc},\ref{diffcbro} holds, it is possible to detect in constant time whether $x_{ij}$ would be strictly positive if a pivoting operation is applied in Phase 1 of \IIO.
\end{Proposition}

\begin{proof}
We recall that, if some degenerate basic edges exist in the spanning tree, they must be in an even position in the unique path from source $i$ to destination $j$ to have $x_{ij} > 0$ with a pivoting operation. Notice also that any path from a source to a destination requires an odd number of edges, while any path between two sources (or two destinations) requires an even number of edges. In case \ref{samec}, checking if $i$ and $j$ have the same color takes constant time. All edges in the path from $i$ to $j$ are non-degenerate, that is all corresponding basic variables have value $>0$. Hence, a pivoting operation would necessarily yield $x_{ij} > 0$. An example of this case in Figure \ref{fig:example} occurs if we consider 
the non-basic variable $x_{42}$.

In case \ref{diffcpc},  there is only one degenerate edge in the path from $i$ to $j$ linking $T_l$ and $T_m$. If subtree $T_l$ is the parent of subtree $T_m$ and the root of $T_m$ is a destination (see, e.g., $T_1$ and $T_2$ in Figure \ref{fig:example} and the non-basic variable $x_{25}$), the degenerate edge is in an odd position in the subpath from source $i$ to the root of  $T_m$ and, therefore, also in the path from $i$ to $j$. Alternatively, if  the root of $T_m$ is a source, the degenerate edge is in an even position in the path from $i$ to $j$. Hence, just checking the root of $T_m$ allows us to evaluate if  $x_{ij} > 0$ would hold with a pivoting operation, and this check can be done in constant time. Similarly, if $T_l$ is a child of $T_m$, condition $x_{ij} > 0$ can hold only if the root of  $T_l$ is a destination. In fact, in this case an odd number of edges connects $i$ to the root of $T_l$ and the next edge is the degenerate edge that is then in an even position in the path from $i$ to $j$. Again, it suffices to check the root of $T_l$ in constant time to evaluate if $x_{ij} > 0$ holds.


 In case \ref{diffcbro}, $T_l$ and $T_m$ are siblings. Let $T_q$ be the parent of $T_l$ and $T_m$ (in Figure \ref{fig:example}, for instance, $T_3$ and $T_4$ are
siblings being children of $T_2$). Correspondingly, there are two degenerate edges in the path from $i$ to $j$. The first degenerate edge links $T_l$ and $T_q$ and the second degenerate edge links $T_q$ and $T_m$. A similar analysis to the previous case implies that both degenerate edges are in an even position in the path from $i$ to $j$ only if the root of $T_l$ is a destination and the root of $T_m$ is a source. Hence, checking  the roots of $T_l$ and $T_m$ in constant time indicates if a pivoting operation would give $x_{ij} > 0$. In the case of Figure \ref{fig:example}, the root of $T_3$ is a source, so, e.g., neither $x_{78} >0$ nor $x_{98}>0$ may occur with a pivoting operation.
\end{proof}

In our algorithmic implementations, as soon as two colors are present, 
we check for pivoting only non-basic variables involving cases \ref{samec},\ref{diffcpc},\ref{diffcbro} while we exclude the non-basic variables of case \ref{diffcboh} to avoid the risk of computing unnecessary cycles. This choice, even if possibly affecting the selection of some promising non-basic variables in the iterations of \IIO, strongly improved the performances of the algorithm.

\section{Computational results}
\label{sec:expe}
In this section, we report the outcome of the computational testing performed to assess the performances of IIO.
The experiments were run as a single thread process on a personal computer equipped with a \textit{11th Gen Intel Core i7-1165G7 2.80GHz} 
processor and 16GB of RAM, and running Ubuntu 20.04.5 LTS.
In a first set of instances, we generated supply/demand quantities and transportation costs from uniform discrete distributions. As mentioned in \cite{Schwinn}, 
transportation instances with uniformly distributed transportation costs are typically among the hardest instances to solve.
In our generation scheme, supply/demand quantities were drawn from the discrete uniform distribution $U(1, 1000)$.
We generated square instances with $K$ sources and $K$ destinations, with $K = 1000, 2000,
4000, 6000, 8000, 10000, 12000, 16000$ and where the transportation costs were drawn from the discrete uniform distribution $U(1, K)$. Ten instances were generated for each value of $K$.
We also considered rectangular instances with $|M| < |N|$ where $|M||N|=36\times10^6$, with $M = 4000, 3000, 2000, 1000$ and $N = 9000, 12000, 18000, 36000$. Also in this case, ten instances were generated for each pair ($|M|,|N|$).


As indicated in the previous section, in applying IIO we consider a parameter $\alpha$ representing the size of the subset of variables with least transport costs for which the reduced costs are computed. After preliminary testing, a reasonable value for instances with uniform distributions was found to be $\alpha=10(|M|+|N|)$.
In Table \ref{tab:iio}, we compare two versions of IIO, denoted IIO+ and IIO-, to our implementation of the standard network simplex algorithm for TP, denoted \NSBDCS, where, at each iteration, the first variable with negative reduced cost is considered for a pivoting operation. IIO+ employs the coloring of the spanning trees associated with basic solutions (as discussed in subsection \ref{ss:coloredspat}) while IIO- does not. We considered the same initial solution, given by the MMR heuristic, and the same value of $\alpha$ in the three algorithms. In Table \ref{tab:iio}, for the two versions of IIO, we list the total number of pivoting operations in both phases (\emph{Pivots}), the total number of macro-iterations (\emph{Macro-it}), i.e., the number of times Phase 1 and Phase 2 are iterated,  and the computational time in seconds (\emph{Time}).  Similarly, for \NSBDCS, we list the total number of pivoting operations (equal to the number of iterations of the algorithm) and the computational time in seconds.
All data are averaged over the ten instances generated for each size.
\begin{table}[h!]
	\scriptsize
	\centering
	\begin{threeparttable}	
		\begin{tabular}{l r r r r r r r r}
			\hline
Instance size		 &  \multicolumn{3}{c}{IIO+} &  \multicolumn{3}{c}{IIO-} & 
\multicolumn{2}{c}{\NSBDCS} \\ \hline
 $|M|x|N|$      &       Pivots &	  Macro-it  &	Time  &            Pivots &	  Macro-it  &	Time   &     Pivots  &	Time   \\ \hline
1000x1000       &       13806 &	131 &	0.046  &              13852 &	105 &	0.127    &           15932  &	0.162 \\ \hline
2000x2000       &      33420	 &	204 &	0.175   &       	   34270 &	170 &	0.543   &            44352  &	1.030 \\ \hline
4000x4000       &      80006	&	328 &	0.699     &          83316 &	266 &	2.399   &           122992  & 7.513 \\ \hline
6000x6000       &      134012 &	468 &	1.812 &             138078 &	388 &	6.508   &           226134 & 29.765 \\ \hline
8000x8000       &      188730 &	558 &	2.912 &             197474 &	454 &	11.683 &            337315  & 61.123 \\ \hline
10000x10000       &      249470 &	672 &	4.360 &             260502 &	551 &	18.716   &          461101  & 106.247 \\ \hline
12000x12000       &      312188 &	773 &	6.571 &             330186 &	627 &	28.034 &            605820  & 178.069 \\ \hline
16000x16000       &      447086 &	946 &	15.760  &           471646 &	783 &	57.227  &           925915  & 405.260 \\ \hline \hline 
4000x9000        &  137000    &	468 & 1.840	  & 141850 & 380	 & 6.723	 &  229983 & 30.346  \\ \hline
3000x12000       &  146660    &503	 & 2.052	  &  151326 & 403	 & 7.223	 &  254207 &  35.016 \\ \hline
2000x18000         & 165720     & 526	 & 2.395	  &  169842 & 431	 & 8.260	 &  299926 &  52.525 \\ \hline
1000x36000         &  216628    & 636	 & 3.973	  &  220548 & 504	 & 13.520	 &   426109 &  128.539 \\ \hline
\end{tabular}
	\caption{Performance comparison of the two versions of IIO and \NSBDCS.}
	\label{tab:iio}
\end{threeparttable}	
\end{table}


From Table \ref{tab:iio}, we evince that IIO- strongly outperforms algorithm \NSBDCS in terms of running times.
IIO+ further improves the performances of IIO-, highlighting the role of the spanning trees coloring. 
IIO+ is more than $10$ times faster than \NSBDCS on instances of size 4000x4000 and larger. 
Interestingly enough, also the number of pivoting operations is inferior in IIO+ compared to \NSBDCS, and the difference between the two values increases as the instances size increases.
We remark that IIO+ manages to solve instances of size 16000x16000 in less than $16$ seconds, and that the behavior of the algorithm is only marginally affected on rectangular instances. The average computational time required by IIO+ on instances of size 1000x36000 is approximately twice the computational time required on instances of size 6000x6000. Finally, we remark that the selection of the initial basic solution does not appear to be too relevant. Additional testing not presented here indicates that, even starting from the well-known NWC method, IIO+ solves instances of size 16000x16000 in approximately 24 seconds on average.


To get a broader overview on the performance of IIO+, we present further details of the results on the square instances in Table \ref{tab:iio_details}. In the table, we report the average length of the paths (\emph{P-length}) computed in Phase 1 and in Phase 2 over the related pivoting operations, the average number of nodes changing color (\emph{Colored nodes}) in Phase 1 over the related pivoting iterations, and the average number of nodes involved in the related pivoting operations (\emph{Involved nodes})  in Phase 2. 

\begin{table}[h!]
	\scriptsize
	\centering
	\begin{threeparttable}	
		\begin{tabular}{l r r r r r}
			\hline
Instance size		 &  $|M|+|N|$  & \multicolumn{4}{c}{IIO+} \\ \hline
       &   &    Phase 1 &	Phase 2  &	Phase 1  &       Phase 2 \\ \hline 
      &   &    P-length  &	P-length   &	Colored nodes  &  Involved nodes  \\ \hline
1000x1000 & 2000 & 31 & 43 & 187 & 49 \\ \hline
2000x2000 & 4000 & 39 & 56 & 291 & 62 \\ \hline
4000x4000 & 8000 & 48 & 73 & 456 & 78 \\ \hline
6000x6000 & 12000 & 55 & 86 & 590 & 91 \\ \hline
8000x8000 & 16000 & 59 & 95 & 709 & 100 \\ \hline
10000x10000 & 20000 & 63 & 103 & 816 & 108 \\ \hline
12000x12000 & 24000 & 67 & 109 & 913 & 114 \\ \hline
16000x16000 & 32000 & 72 & 121 & 1091 & 125 \\ \hline
\end{tabular}
	\caption{Additional results of IIO+ algorithm.}
	\label{tab:iio_details}
\end{threeparttable}	
\end{table}

We notice that all entries in the columns 3-6 of Table \ref{tab:iio_details} are from one to two orders of magnitude inferior to $|M|+|N|$. Correspondingly, in practice the average complexity of a pivoting operation of IIO+ is much lower than the complexity of processing a pivoting operation of \NSBDCS, which necessarily requires the computation of all $u_i$ and $v_j$ multipliers in $\Theta(|M|+|N|)$.

We compared then IIO+ on the same instances with various benchmark solvers available from the literature.
In particular, we considered the solver Cplex (version  20.1) applying the primal simplex, the dual simplex and the network simplex. Similarly, we tested the solver Gurobi (version  10.0) applying the primal simplex, the dual simplex and the network simplex. Preliminary testing indicated that the barrier method both for Cplex and Gurobi was much less efficient than the other methods, so we did not consider this further method in the performance comparison. Also, the best performance of both solvers were obtained by deactivating the presolving option and activating the sifting option. The solvers ran with a single thread while the other parameters were set to their default values. 
Besides, we considered an efficient C++ implementation of the network simplex developed in \cite{Bonneel}, hereafter denoted \NSB, which is based on the graph library \cite{LEMON}. Since the behavior of  IIO+ on rectangular instances is similar to the one on square instances, we limited the analysis to square instances. 
For all solution methods, we report in Table \ref{tab:iiovssolvers} the average computational time over 10 instances for each category.

\begin{table}[h!]
	\scriptsize
	\centering
	\begin{threeparttable}	
		\begin{tabular}{l r r r r r r r r }
			\hline
 Instance size	&  \multicolumn{3}{c}{Cplex} &  \multicolumn{3}{c}{Gurobi} & 
\NSB  & IIO+  \\ \hline

       &       Primal &	 Dual  &	Network  &           Primal &	 Dual  &	Network  &           	&   \\ \hline
       &       Time &	 Time &	 Time &	 Time &	 Time &	 Time &	 Time &	Time \\ \hline
1000x1000 & 1.670 & 2.433 & 1.466 & 2.020 & 0.738 & 2.163 & 0.252 & 0.046 \\ \hline
2000x2000 & 7.546 & 16.351 & 6.923 & 188.473 & 3.245 & 199.772 & 0.881 & 0.175 \\ \hline
4000x4000 & 34.692 & 122.782 & 35.858 & 48.259 & 13.281 & 53.531 & 3.677 & 0.699 \\ \hline
6000x6000 &  85.775 & 280.27  & 105.642 & 254.598 & 30.636 & 290.219 & 9.985 & 1.812 \\ \hline
\hline 
8000x8000 & - & - & - & - & - & - & 19.103 & 2.912 \\ \hline
10000x10000 & - & - & - & - & - & - & 33.458 & 4.360 \\ \hline
12000x12000 & - & - & - & - & - & - & 50.336 & 6.571 \\ \hline
16000x16000 & - & - & - & - & - & - & - & 15.760 \\ \hline
\hline 
\end{tabular}
	\caption{Performance comparison of IIO+ and benchmark solvers.}
	\label{tab:iiovssolvers}
\end{threeparttable}	
\end{table}

\normalsize

The results in Table \ref{tab:iiovssolvers} indicate that IIO+ strongly dominates all competitors among which \NSB shows the best performance. We also note that the larger the instances, the larger the difference in running times between IIO+ and the other solvers. We remark that Cplex and Gurobi were limited to a size not superior to  6000x6000 due to memory requirements. The same consideration holds for \NSB and instances 12000x12000. Further, while the dual simplex of Gurobi shows up to be quite efficient among the competitors, we observe an unexpected poor performance of the primal simplex and the network simplex of Gurobi on instances 2000x2000.


We also considered the square benchmark instances of TP introduced in \cite{Schrieber}, denoted as DOTmark instances. As mentioned in Section \ref{Intro}, these instances derive from applications in image processing, and their transportation costs are determined according to specific geometrical structures. In particular, sources and destinations are located on grids with corresponding coordinates, and each unit transportation cost $c_{ij}$ between a source $i$ and a destination $j$ is computed as the squared Euclidean distance between $i$ and $j$. We refer to \cite{Schrieber} for further details on the matter.
Among the methods analyzed in \cite{Schrieber} for TP, the best performances in terms of computational times were obtained by the AHA method \citep{AHA, Merigot}, here denoted as AHAM, which is actually not an exact method. The best exact approach, tailored for this type of instances, was the so-called Shielding Neighborhood Method (here denoted as SNM) presented in \cite{Schmitzer}, where restricted instances of the original problem are iteratively solved to optimality until the last solution can be shown to be optimal for the initial problem. Each restricted instance contains a subset of the variables of the initial problem constituted by a  feasible solution and a peculiar (typically small) neighborhood of that solution. The approach iteratively moves from a feasible solution to a better one until a solution is optimal for two consecutive iterations, implying that the computed solution is also optimal to the original problem. In \cite{Schrieber}, each restricted instance was solved with the network simplex of Cplex. 
We tested the performances of IIO+ by applying it as a subroutine of the Shielding method to solve the sequences of restricted transportation instances (instead of using the network simplex of Cplex). For each restricted instance, we set the value of $\alpha$ equal to the size of the instance. Also, we defined a threshold limit $\beta$ on the number of macro-iterations of IIO+. After preliminary testing, we set $\beta = 28$. This choice led to some savings in computational times in the search for a global optimal solution. However, we remark that the application of IIO+ without $\beta$ also showed very promising results. 
We denoted the corresponding approach as SNM with IIO+. We also tested \NSB that specifically exploits the structure of the transportation costs on instances such as the DOTmark instances. We report the corresponding results in Table \ref{tab:shielding}. 
The tests in \cite{Schrieber} were performed on a Linux server (AMD Opteron Processor 6140 from 2011 with 2.6 GHz). 
We report here the computation times (taken from \cite{Schrieber}) of SNM and AHAM 
on that machine that is comparable to the machine used in our tests. 
An asterisk in Table \ref{tab:shielding} points out that the related times refer to a different machine.


\begin{table}[H]
	\scriptsize
	\centering
	\begin{threeparttable}	
		\begin{tabular}{ l r r r r r r r}
			\hline
Instance type	&	 Instance size		 &  \multicolumn{3}{c}{SNM with IIO+} &  \NSB & AHAM & SNM \\ \hline
   &    &       Pivots &	Macro-it  &	Time  & Time  & Time$^*$ & Time$^*$ \\ \hline
WhiteNoise-32       & 1024x1024 &   31810     &	371 & 0.101  &  0.233 & 3.28 & 0.67 \\ \hline
GRFrough-32       & 1024x1024 &    41702    &  447	&0.120  &  0.222  & 3.19 & 1.08 \\ \hline
GRFmoderate-32  & 1024x1024      &   41437     & 524	& 0.129 &   0.245 & 3.17 & 1.86 \\ \hline
GRFsmooth-32     & 1024x1024   &    44298    & 597	& 0.139 &    0.259 & 4.39 & 2.66\\ \hline
LogGRF-32     & 1024x1024  &   52648     &817	& 0.169 &  0.273  & 6.80 & 3.00\\ \hline
LogitGRF-32   & 1024x1024    &   48833     &721	& 0.158 &   0.269 & 8.49 & 2.40\\ \hline
CauchyDensity-32 & 1024x1024      & 36927  & 469 & 0.113 & 0.250    & 3.76 & 3.62\\ \hline
Shapes-32   & 1024x1024   &    16192    &	272 & 0.042 &  0.077  & 1.27 & 0.92\\ \hline
ClassicImages-32 & 1024x1024         &  41277      & 539	& 0.136 & 0.254    & 2.01 & 1.58\\ \hline
Microscopy-32    & 1024x1024    &   35230     &562	& 0.152  &   0.137 & 3.14 & 1.66\\ \hline \hline \hline
Overall-32      & 1024x1024    &   39035      & 532	& 0.126  &   0.222 & 3.95 & 1.94\\ \hline\hline \hline
WhiteNoise-64       & 4096x4096 &   305972     & 1286	& 1.475 &  4.064  & 36 & 13\\ \hline
GRFrough-64       & 4096x4096 &   356909     & 1537	&  1.798 &  10.696  & 20 & 24\\ \hline
GRFmoderate-64  & 4096x4096      &   420092     & 2461	&  2.248&   6.432 & 23 & 51\\ \hline
GRFsmooth-64     & 4096x4096   &     558740   & 4130	& 3.361 &   7.598 & 57 & 80\\ \hline
LogGRF-64     & 4096x4096   &    695497    &  5256	&  4.158 &   8.252 & 59 & 79\\ \hline
LogitGRF-64   & 4096x4096     &   489050     & 2959	& 2.564 & 7.321   & 77 & 69\\ \hline
CauchyDensity-64 & 4096x4096       &   439535 	&  2410   & 2.213  & 5.386    & 30 & 97\\ \hline
Shapes-64  & 4096x4096    &    152706    & 1023	&  0.064 &   1.159 & 12 & 25\\ \hline
ClassicImages-64 & 4096x4096          &   393515     & 2438	&2.184  &  6.406   & 18 & 41\\ \hline
Microscopy-64    & 4096x4096     &    89720    & 795	&  3.487 &   0.613 & 24 & 26\\ \hline \hline \hline
Overall-64      & 4096x4096     &   390174     &2429	&   2.413 &   5.793 & 36 & 51\\ \hline \hline \hline
\end{tabular}
	\caption{Comparing SNM with IIO+ to \NSB, AHAM and SNM.}
	\label{tab:shielding}
\end{threeparttable}	
\end{table}

\normalsize

We notice that \NSB performs much better than AHAM and SNM regardless the 
marginal differences among the machines used. Still, Shielding with IIO+ is, on the average, 
approximately twice faster than \NSB.

\section{Conclusions}
\label{sec:concl}
We have proposed \IIO, a new exact algorithm for the transportation problem.
The algorithm requires in input a basic feasible solution and is composed by two main phases that are iteratively repeated until an optimal basic solution is computed. In the first phase, the algorithm progressively improves the current solution by moving through pivoting operations towards a feasible solution interior to the constraint set polytope. The second phase moves back towards a further improved basic solution always by means of pivoting operations. 
An extensive computational campaign showed the efficiency of \IIO, which strongly outperforms all the current state-of-the-art approaches.
Several issues may be worthy to investigate in future research.
First, the proposed approach can be extended to solve the minimum cost flow problem. While a natural extension of \IIO is easily conceivable, an efficient implementation would require a detailed analysis of the features of the minimum cost flow problem.
Second, it could be worthy to tune the algorithm for solving the assignment problem where the basic solutions are strongly degenerate.
Besides, in this work no special effort has been dedicated to the search of an initial solution
(possibly non-basic, given the availability of Phase 2) for the transportation problem. We could explore new heuristic methods to further reduce the computational time required by \IIO and to quickly compute effective feasible solutions for very large instances of the transportation problem.
Finally, it could be interesting to further study the interplay between \IIO and specific cost functions of transportation instances associated with practical applications in computer vision and machine learning.


\begin{thebibliography}{}

\bibitem[{Aurenhammer et~al.(1998)}]{AHA} 
Aurenhammer F, Hoffmann F, Aronov B (1998) Minkowski-type theorems and leastsquares clustering. {\it Algorithmica} 20(1):61--76.

\bibitem[{Armstrong and Jin (1997)}]{Armstrong} 
Armstrong RD, Jin Z (1997) A new strongly polynomial dual network simplex algorithm. {\it Mathematical Programming} 78:131--148.

\bibitem[{Bassetti et~al.(2020)}]{BaGuaVe} 
Bassetti F, Gualandi S, Veneroni M (2020) On the Computation of Kantorovich--Wasserstein Distances Between Two-Dimensional Histograms by Uncapacitated Minimum Cost Flows. {\it SIAM Journal on Optimization} 30(3):2441--2469.

\bibitem[{Bonneel et~al.(2011)}]{Bonneel} 
Bonneel N, van de Panne M, Paris S, Heidrich W (2011) Displacement interpolation using lagrangian mass transport. {\it ACM Transactions on Graphics} 30(6):1--12.

\bibitem[{Bulut(2017)}]{Bulut} 
Bulut H (2017) Multiloop transportation simplex algorithm. {\it Optimization Methods and Software} 32(6):1206--1217.

\bibitem[{Cunningham(1976)}]{Cunningham} 
Cunningham W (1976) A network simplex method.
{\it Mathematical Programming} 11:105--116.

\bibitem[{Dantzig(1951)}]{Dantzig} 
Dantzig GB (1951) Application of the simplex method to a transportation problem. In Koopmans TC, ed., {\it Activity
 Analysis of Production and Allocation} 13:359--373, John Wiley and Sons.

\bibitem[{Ford and Fulkerson(1956)}]{FordFulkerson} 
Ford LR, Fulkerson DR (1956) Solving the transportation problem. {\it Management Science} 3(1):24--32.

\bibitem[{Gottschlich and Schuhmacher(2014)}]{Gottschlich} Gottschlich C, Schuhmacher D (2014) The shortlist method for fast computation of the earth mover's distance and finding optimal solutions to transportation problems. {\it PLoS ONE} 
9(10):e110214.

\bibitem[{Hitchcock(1941)}]{Hitchcock} 
Hitchcock FL (1941) The Distribution of a Product from Several Sources to Numerous Localities. {\it Journal of Mathematics and Physics} 20:224--230.

\bibitem[{Kantorovich(1942)}]{Kantorovich} 
Kantorovich LV (1942) On the translocation of masses. {\it C.R. (Doklady) Acad. Sci. URSS (N.S.)} 37:199--201.

\bibitem[{Kov\'acs(2015)}]{Kovacs} 
Kov\'acs P (2015) Minimum-cost flow algorithms: an experimental evaluation. {\it Optimization Methods and Software} 30:94--127. 

\bibitem[{LEMON(2010)}]{LEMON} (2010) LEMON. Library for efficient modeling and optimization in networks. http://lemon.cs.elte.hu/trac/lemon.

\bibitem[{Luenberger and Ye(2008)}]{Luenberger} 
Luenberger DG, Ye Y (2008) Linear and Nonlinear programming. {\it International Series in Operations Research and Management Science}, Springer, 3rd Edition.

\bibitem[{M\'erigot(2011)}]{Merigot} 
M\'erigot Q (2011) A multiscale approach to optimal transport. {\it Computer Graphics Forum} 30(5):1583--1592.

\bibitem[{Monge(1781)}]{Monge} 
Monge G (1781) M\'emoire sur la th\'eorie des  d\'eblais et des remblais. In {\it Histoire de l'Acad\'emie Royale des Sciences de Paris} 666--704.

\bibitem[{Orlin(1996)}]{Orlin} 
Orlin JB (1996) A polynomial time primal network simplex algorithm for minimum cost flows. {\it Mathematical Programming} 78(2):109--129.


\bibitem[{Schmitzer(2016)}]{Schmitzer} 
Schmitzer B (2016) A sparse multiscale algorithm for dense optimal transport. {\it Journal of  Mathematical Imaging and Vision} 56(2):238--259.

\bibitem[{Schrieber et~al.(2016)}]{Schrieber} 
Schrieber J, Schuhmacher D, Gottschlich C (2016) DOTmark - a benchmark for discrete optimal transport. {\it IEEE Access} 5:271--282.

\bibitem[{Schwinn and Werner(2019)}]{Schwinn}
Schwinn J, Werner R (2019) On the effectiveness of primal and dual heuristics for the transportation problem. {\it IMA Journal of Management Mathematics} 30:281--303.
\end{thebibliography}
\end{document}